\documentclass[11pt, twoside, eqno]{article}
\usepackage{latexsym}
\usepackage{amssymb}
\usepackage{amsfonts}
\begin{document}
\begin{center}
{\bf \Large On Fuzzy Ideals and Level Subsets of
Ordered $\Gamma$-Groupoids}

\medskip

\noindent{\bf Niovi Kehayopulu}\\
{\it Department of Mathematics, University of Athens,
15784 Panepistimiopolis, Greece}\\email: nkehayop@math.uoa.gr

\noindent October 15, 2014\end{center}
\bigskip{\small

\noindent{\bf Abstract.} We characterize the fuzzy left (resp. right) 
ideals, the fuzzy ideals and the fuzzy prime (resp. semiprime) ideals 
of an ordered $\Gamma$-groupoid $M$ in terms of level subsets and we 
prove that the cartesian product of two fuzzy left (resp. right) 
ideals of $M$ is a fuzzy left (resp. right) ideal of $M\times M$, and 
the cartesian product of two fuzzy prime (resp. semiprime) ideals of 
$M$ is a fuzzy prime (resp. semiprime) ideal of $M\times M$. As a 
result, if $\mu$ and $\sigma$ are fuzzy left (resp. right) ideals, 
ideals, fuzzy prime or fuzzy semiprime ideals of $M$, then the 
nonempty level subsets $(\mu\times\sigma)_t$ are so.\bigskip

\smallskip

 }
\section{Introduction and prerequisites}For two nonempty sets $M$ and 
$\Gamma$, we denote by the letters of the English alphabet the 
elements of $M$ and by the letters of the Greek alphabet the elements 
of $\Gamma$, and define $M\Gamma M:=\{a\gamma b \mid a,b\in M, 
\gamma\in\Gamma\}$. Then $M$ is called a {\it $\Gamma$-groupoid} if

(1) $M\Gamma M\subseteq M$ and

(2) if $a,b,c,d\in M$ and $\gamma,\mu\in\Gamma$ such that $a=c$, 
$b=d$ and $\gamma=\mu$, then

\hspace{0.5cm} have $a\gamma b=c\mu d$.\\If, in addition, for all 
$a,b\in M$ and all $\gamma,\mu\in\Gamma$, $a\gamma (b\mu c)=(a\gamma 
b)\mu c$, then $M$ is called a {\it $\Gamma$-semigroup}. For 
$\Gamma=\{\gamma\}$, the $\Gamma$-groupoid $M$ is the so called 
groupoid and the ordered $\Gamma$-groupoid the ordered groupoid (: 
$po$-groupoid), the $\Gamma$-semigroup is the semigroup and the 
ordered $\Gamma$-semigroup the ordered semigroup (where $\gamma$ is 
the multiplication on $M$ usually denoted by ``$\cdot$").
An {\it ordered $\Gamma$-groupoid} (: {\it $po$-$\Gamma$-groupoid}) 
is a $\Gamma$-groupoid $M$ with an order relation ``$\le$" on $M$ 
such that $a\le b$ implies $a\gamma c\le b\gamma c$ and $c\gamma a\le 
c\gamma b$ for all $\gamma\in \Gamma$. Let $M$ be a 
$po$-$\Gamma$-groupoid. A subset $T$ of $M$ is called {\it prime} if 
for every $a,b\in M$ and every $\gamma\in\Gamma$ such that $a\gamma 
b\in T$, we have $a\in T$ or $b\in T$. A subset $T$ of $M$ is called 
{\it semiprime} if for every $a\in M$ and every $\gamma\in\Gamma$ 
such that $a\gamma a\in T$, we have $a\in T$. A nonempty subset $A$ 
of $M$ is called a {\it left} (resp. {\it right}) {\it ideal} of $M$ 
if

(1) $M\Gamma A\subseteq A$ (resp. $A\Gamma M\subseteq A$) and

(2) if $a\in A$ and $M\ni b\le a$ implies $b\in A$.\\
If the set $A$ is both a left and a right ideal of $M$, then it is 
called an {\it ideal} of $M$. If $M$ is an ordered $\Gamma$-groupoid, 
then any mapping $\mu : M\rightarrow [0,1]$ is called a {\it fuzzy 
subset} of $M$ (or a {\it fuzzy set} in $M$) (L. Zadeh). The mapping 
$\mu$ is called a {\it fuzzy left ideal} of $M$ if

(1) $\mu(x\gamma y)\ge \mu(y)$ for every $x,y\in M$ and every 
$\gamma\in\Gamma$ and

(2) if $x\le y$ implies $\mu(x)\ge \mu (y)$.\\It is called a {\it 
fuzzy right ideal} of $M$ if

(1) $\mu(x\gamma y)\ge \mu(x)$ for every $x,y\in M$ and every 
$\gamma\in\Gamma$ and

(2) if $x\le y$ implies $\mu(x)\ge \mu (y)$.\\A fuzzy subset which is 
both a fuzzy left and a fuzzy right ideal of $M$ is called a {\it 
fuzzy ideal} of $M$. A fuzzy subset $\mu$ of $M$ is a fuzzy ideal of 
$M$ if and only if

(1) $\mu(x\gamma y)\ge\max\{\mu(x),\mu(y)\}$ for every $x,y\in M$ and 
every $\gamma\in\Gamma$ and

(2) if $x\le y$ implies $\mu(x)\ge \mu (y)$.

A characterization of fuzzy prime subsets of a semigroup in terms of 
level subsets has been considered in [4; Lemma 2.3] in which 
$\lambda$ should be replaced by $t$ (or $t$ should be replaced by 
$\lambda$) and the commutativity of the semigroup is not necessary as 
the same holds in semigroups is general. As an immediate consequence 
of the Lemma 2.3 in [4], a fuzzy ideal of a semigroup $S$ is prime if 
and only if for every $t\in [0,1]$ the $t$-level subset $f_t:=\{x\in 
S \mid f(x)\ge t\}$ (of $S$), if it is nonempty, is a prime ideal of 
$S$. This is the Theorem 3.1 in [3] in which the ``$xy\subseteq 
\mu_t$" (in the second line of the proof) should be replaced by 
``$xy\in \mu_t$" and the proof of the converse statement (lines 5--11 
of the proof) should be corrected. A characterization of fuzzy 
semiprime ideals of a semigroup in terms of level subsets has been 
considered in the Theorem 3.2 in [3] but the proof of the ``converse" 
statement in it should be corrected. The cartesian product of two 
fuzzy left (resp. fuzzy right) ideals of semigroups and the cartesian 
product of two fuzzy prime (resp. fuzzy semiprime) ideals of a 
semigroup has been studied in [3]. On the other hand, a 
characterization of fuzzy prime and fuzzy semiprime ideals of ordered 
semigroups in terms of level subsets has been considered in the 
Theorems 2.6 and 2.7 in [6] from which the Theorems 3.1 and 3.2 in 
[3] are also obtained. For a characterization of fuzzy ideals of 
ordered semigroup in terms of level subsets see the Lemma 2.4 in [6] 
and the Lemma 2.7 in [5]. The reference in Lemma 2.4 in [6] should be 
corrected. In the ``converse statement" of the proof of Theorem 3.1 
in [3] as well as in the proof of Theorem 3.2 in [3], the phrase 
``Let every nonempty subset $\mu_t$ of $\mu$ be a prime (semiprime) 
ideal of $S$" is better to be replaced by the phrase ``Let every 
subset $\mu_t$ of $\mu$ be a prime (semiprime) ideal of $S$", and 
this is because the ideals are, by definition, nonempty sets. 
Finally, the proofs of Propositions 4.4 and 4.6 in [3] should be 
omitted as they are immediate consequences of Propositions 4.2 and 
4.3 and the Theorem 3.2 given in the same paper.

In the present paper we first characterize the fuzzy left (right) 
ideals, the fuzzy ideals, the fuzzy prime and the fuzzy semiprime 
ideals of an ordered $\Gamma$-groupoid in terms of level subsets. 
Then we prove that the cartesian product of two fuzzy left (resp. 
fuzzy right) ideals of an ordered $\Gamma$-groupoid $M$ is a fuzzy 
left (resp. fuzzy right) ideal of $M\times M$. Thus the cartesian 
product of two
fuzzy ideals of $M$ is a fuzzy ideal of $M\times M$. Moreover, the 
cartesian product of two
fuzzy prime (resp. fuzzy semiprime) ideals of a $\Gamma$-groupoid $M$ 
is a fuzzy prime (resp. fuzzy semiprime) ideal of $M\times M$. As a 
consequence, if $\mu$ and $\sigma$ are fuzzy left (resp. fuzzy right) 
ideals of an ordered $\Gamma$-groupoid $M$ then, for any $t\in[0,1]$ 
if the level subset $(\mu\times\sigma)_t$ is nonempty, then it is a 
left (resp. right) ideal of $M\times M$. If $\mu$ and $\sigma$ are 
fuzzy ideals of $M$ and the level subset $(\mu\times\sigma)_t$ is 
nonempty, then it is an ideal of $M\times M$. If $\mu$ and $\sigma$ 
are fuzzy prime (resp. fuzzy semiprime) ideals of $M$, then the 
nonempty level subsets $(\mu\times\sigma)_t$  are prime (resp. 
semiprime) ideals of $M\times M$. The present paper serves as an 
example to show the way we pass from fuzzy ordered groupoids (resp. 
fuzzy ordered semigroups) to fuzzy ordered $\Gamma$-groupoids (resp. 
fuzzy ordered $\Gamma$-semigroups) and from fuzzy groupoids (resp. 
fuzzy semigroups) to fuzzy $\Gamma$-groupoids (resp. fuzzy 
$\Gamma$-semigroups). On the other hand, from the results of 
$\Gamma$-groupoids or $\Gamma$-ordered groupoids where 
$\Gamma=\{\gamma\}$ ($\gamma$ being a symbol) the corresponding 
results on groupoids or ordered groupoids are obtained. The fuzzy 
sets in ordered groupoids have been introduced in [1] and one can 
find several papers on fuzzy ordered semigroups in the bibliography.
\section{Characterization of prime and semiprime fuzzy ideals in 
terms of level subsets}Following the terminology of fuzzy prime 
subset of a groupoid introduced in [1, 2], we give the following 
definition\medskip

\noindent{\bf Definition 1.} Let $M$ be an ordered $\Gamma$-groupoid 
(or a $\Gamma$-groupoid). A fuzzy subset $\mu$ of $M$ is called {\it 
fuzzy prime subset of $M$} or {\it prime fuzzy subset of $M$} 
if$$\mu(x\gamma y)\le\max\{\mu(x),\mu(y)\}$$for all $x,y\in M$ and 
all $\gamma\in\Gamma$.\\Recall that if $\mu$ is a fuzzy prime ideal 
of $M$, then for every $x,y\in M$ and every $\gamma\in\Gamma$, we 
have $\mu(x\gamma y)=\max\{\mu(x),\mu(y)\}$. So a fuzzy ideal $\mu$ 
of $M$ can be called {\it prime} if $\mu(x\gamma 
xy)=\max\{\mu(x),\mu(y)\}$ for all $x,y\in M$ and all 
$\gamma\in\Gamma$.\medskip

\noindent{\bf Definition 2.} Let $M$ be an ordered $\Gamma$-groupoid 
(or a $\Gamma$-groupoid). A fuzzy subset $\mu$ of $M$ is called {\it 
fuzzy semiprime subset of $M$} or {\it semiprime fuzzy subset of $M$} 
if$$\mu(x)\ge \mu(x\gamma x)$$for every $x\in M$ and every 
$\gamma\in\Gamma$.\medskip

\noindent{\bf Notation 3.} If $\mu$ is a fuzzy subset of an ordered 
$\Gamma$-groupoid (or a $\Gamma$-groupoid) $M$ then, for any $t\in 
[0,1]$ (: the closed interval of real numbers), we denote by $\mu_t$ 
the subset of $M$ defined by$$\mu_t:=\{x\in M \mid \mu(x)\ge t\}.$$ 
The set $\mu_t$ is called the $t$-level subset or just level subset 
of $\mu$.\medskip

\noindent{\bf Theorem 4.} {\it Let M be an ordered $\Gamma$-groupoid. 
If $\mu$ is a fuzzy left ideal of M and $\mu_t\not=\emptyset$, then 
$\mu_t$ is a left ideal of $M$. ``Conversely", if $\mu_t$ is a left 
ideal of $M$ for every $t$, then $\mu$ is a fuzzy left ideal of 
M}.\medskip

\noindent{\bf Proof.} $\Longrightarrow$. Suppose $\mu$ is a fuzzy 
left ideal of $M$ and $\mu_t\not=\emptyset$ for some $t\in[0,1]$. 
Then $M\Gamma\mu_t\subseteq \mu_t$. Indeed: Let $a\in M$, 
$\gamma\in\Gamma$ and $b\in\mu_t$. Since $\mu$ is a fuzzy left ideal 
of $M$, we have $\mu(a\gamma b)\ge \mu(b)$. Since $b\in\mu_t$, we 
have $\mu(b)\ge t$. Then $\mu(a\gamma b\ge t$, and $a\gamma b\in 
\mu_t$.
Let $a\in \mu_t$ and $M\ni b\le a$. Then $b\in \mu_t$. Indeed: Since 
$a\in \mu_t$, we have $\mu(a)\ge t$. Since $b\le a$ and $\mu$ is a 
fuzzy left ideal of $M$, we have $\mu(b)\ge \mu(a)$. Then $\mu(b)\ge 
t$, and $b\in\mu_t$. Thus $\mu_t$ is a left ideal of $M$.\smallskip

\noindent$\Longleftarrow$. Suppose $\mu_t$ is a left ideal of $M$ for 
every $t$ and let $a,b\in M$ and $\gamma\in\Gamma$. Then $\mu(a\gamma 
b)\ge \mu (b)$. Indeed: Since $\mu(b)\in [0,1]$ and $\mu(b)\ge 
\mu(b)$, we have $b\in \mu_{\mu(b)}$. Since $\mu_{\mu(b)}$ is a left 
ideal of $M$, we have $a\gamma b\in M\Gamma \mu_{\mu(b)}\subseteq 
\mu_{\mu(b)}$. Then $a\gamma b\in \mu_{\mu(b)}$, and $\mu(a\gamma 
b)\ge \mu(b)$. Let now $a\le b$. Then $\mu(a)\ge \mu(b)$. Indeed: 
Since $b\in\mu_{\mu(b)}$, $M\ni a\le b$ and $\mu_{\mu(b)}$ is a left 
ideal of $M$, we have $a\in\mu_{\mu(b)}$, then $\mu(a)\ge 
\mu(b)$.$\hfill\Box$\\In a similar way we have the following\medskip

\noindent{\bf Theorem 5.} {\it Let M be an ordered $\Gamma$-groupoid. 
If $\mu$ is a fuzzy right ideal of M and $\mu_t\not=\emptyset$, then 
$\mu_t$ is a right ideal of $M$. ``Conversely", if $\mu_t$ is a right 
ideal of $M$ for every $t$, then $\mu$ is a fuzzy right ideal of 
M}.\\By Theorems 4 and 5, we have the following theorem\medskip

\noindent{\bf Theorem 6.} {\it Let M be an ordered $\Gamma$-groupoid. 
If $\mu$ is a fuzzy ideal of M and $\mu_t\not=\emptyset$ for some 
$t\in[0,1]$, then $\mu_t$ is an ideal of $M$. ``Conversely", if 
$\mu_t$ is an ideal of $M$ for every $t\in [0,1]$, then $\mu$ is a 
fuzzy ideal of M}.\medskip

\noindent{\bf Lemma 7.} {\it Let M be an ordered $\Gamma$-groupoid.
Then $\mu$ is a fuzzy prime subset of $M$ if and only if the level 
subset $\mu_t$ is a prime subset of $M$ for every $t$.}\medskip

\noindent{\bf Proof.} $\Longrightarrow$. Let $a,b\in M$ and 
$\gamma\in\Gamma$ such that $a\gamma b\in\mu_t$. Then $a\in\mu_t$ or 
$b\in\mu_t$. Indeed: Since $a\gamma b\in\mu_t$, we have $\mu(a\gamma 
b)\ge t$. Since $\mu$ is a fuzzy prime subset of $M$, we have 
$\mu(a\gamma b)\le\max\{\mu(a),\mu(b)\}$. Since $\mu(a),\mu(b)\in 
[0,1]$, we have $\mu(a)\le\mu(b)$ or $\mu(b)\le\mu(a)$. If $\mu(a)\le 
\mu(b)$, then $\max\{\mu(a),\mu(b)\}=\mu(b)$, and $t\le\mu(b)$, so 
$b\in\mu_t$. If $\mu(b)\le \mu(a)$, then $t\le \mu(a\gamma 
b)=\mu(a)$, and $a\in\mu_t$.\smallskip

\noindent$\Longleftarrow$. Suppose $\mu_t$ is a prime subset of $M$ 
for every $t$ and let $x,y\in M$ and $\gamma\in\Gamma$. Then 
$\mu(x\gamma y)=\max\{\mu(x),\mu(y)\}$. Indeed: Since $x\gamma 
y\in\mu_{\mu(x\gamma y)}$, by hypothesis, we have 
$x\in\mu_{\mu(x\gamma y)}$ or $y\in\mu_{\mu(x\gamma y)}$. Then 
$\mu(x)\ge\mu(x\gamma y)$ or $\mu(y)\ge\mu(x\gamma y)$, thus
$\max\{\mu(x),\mu(y)\}\ge \mu(x\gamma y)$. $\hfill\Box$\\By Theorem 6 
and Lemma 7, we have the following theorem\medskip

\noindent{\bf Theorem 8.} {\it Let M be an ordered $\Gamma$-groupoid. 
If $\mu$ is a fuzzy prime ideal of M and $\mu_t\not=\emptyset$, then 
$\mu_t$ is a prime ideal of $M$. ``Conversely", if $\mu_t$ is a prime 
ideal of $M$ for every $t$, then $\mu$ is a fuzzy prime ideal of 
M}.\medskip

\noindent{\bf Lemma 9.} {\it Let M be an ordered $\Gamma$-groupoid. 
Then $\mu$ is a fuzzy semiprime subset of M if and only if the level 
subset $\mu_t$ is a semiprime subset of $M$ for every $t$}.\medskip

\noindent{\bf Proof.} $\Longrightarrow$. Let $t\in [0,1]$, $a\in M$ 
and $\gamma\in\Gamma$ such that $a\gamma a\in\mu_t$. Then 
$a\in\mu_t$. Indeed: Since $\mu$ is a fuzzy semiprime subset of $M$, 
we have $\mu(a)\ge \mu(a\gamma a)$. Since $a\gamma a\in\mu_t$, we 
have $\mu(a\gamma a)\ge t$. Then $\mu(a)\ge t$, and 
$a\in\mu_t$.\smallskip

\noindent$\Longleftarrow$. Let $a\in M$ and $\gamma\in\Gamma$. Then 
$\mu(a)\ge \mu(a\gamma a)$. Indeed: By hypothesis, $\mu_{\mu(a\gamma 
a)}$ is a semiprime subset of $M$. Since $a\gamma 
a\in\mu_{\mu(a\gamma a)}$, we have $a\in \mu_{\mu(a\gamma a)}$, then 
$\mu(a)\ge\mu(a\gamma a)$, so $\mu$ is fuzzy 
semiprime.$\hfill\Box$\\By Theorem 6 and Lemma 9, we have the 
following theorem\medskip

\noindent{\bf Theorem 10.} {\it Let M be an ordered 
$\Gamma$-groupoid. If $\mu$ is a fuzzy semiprime ideal of M and 
$\mu_t\not=\emptyset$, then $\mu_t$ is a semiprime ideal of $M$. 
``Conversely", if $\mu_t$ is a semiprime ideal of $M$ for every $t$, 
then $\mu$ is a fuzzy semiprime ideal of M}.\medskip

As a consequence, given a groupoid or an ordered groupoid $G$, if 
$\mu$ is a fuzzy left ideal, fuzzy right ideal, fuzzy ideal, fuzzy 
prime ideal or fuzzy semiprime ideal of $G$, respectively, then the 
nonempty level subsets $\mu_t$ of $\mu$ are left ideals, right 
ideals, ideals, prime ideals or semiprime ideals of $G$, 
respectively. ``Conversely" if for a fuzzy subset $\mu$ of $G$ the 
and any $t\in [0,1]$ the level subset $\mu_t$ of
$\mu$ is a left ideal, right ideal, ideal, prime ideal or semiprime 
ideal of $G$, respectively, then $\mu$ is a fuzzy left ideal, fuzzy 
right ideal, fuzzy ideal, fuzzy prime ideal or fuzzy semiprime ideal 
of $G$, respectively.
\section{Cartesian product of fuzzy ideals, fuzzy prime and fuzzy 
semiprime ideals}If $(M,\le,\Gamma)$ is an ordered $\Gamma$-groupoid, 
$M\times M:=\{(x,y) \mid x,y\in M\}$ and for any $(a,b),(c,d)\in 
M\times M$ and any $\gamma\in\Gamma$ we define $$(a,b)\gamma 
(c,d):=(a\gamma c,b\gamma d),$$ then $(M\times M,\le,\Gamma)$ is an 
ordered $\Gamma$-groupoid as well.

For two fuzzy subsets $\mu$ and $\sigma$ of an ordered 
$\Gamma$-groupoid $M$, the cartesian product of $\mu$ and $\sigma$ is 
the fuzzy subset of $M\times M$ defined by
$$\mu\times \sigma : M\times M \rightarrow [0,1] \mid 
(x,y)\rightarrow \min\{\mu(x),\sigma(y)\}.$$That is,$$(\mu\times 
\sigma){\Big(}(x,y){\Big)}:=\min\{\mu(x),\sigma(y)\}$$for every 
$x,y\in M$. As no confusion is possible, we write $(\mu\times 
\sigma)(x,y)$ instead of $(\mu\times 
\sigma){\Big(}(x,y){\Big)}$.\medskip

\noindent{\bf Theorem 11.} {\it Let M be an ordered $\Gamma$-groupoid 
and $\mu$, $\sigma$ fuzzy left (resp. right) ideals of M. Then 
$\mu\times\sigma$ is a fuzzy left (resp. right) ideal of $M\times 
M$.}\medskip

\noindent{\bf Proof.} Suppose $\mu$ and $\sigma$ are fuzzy left 
ideals of $M$. Let $(a,b),(c,d)\in M\times M$ and $\gamma\in\Gamma$. 
Then $$(\mu\times \sigma){\Big(}(a,b)\gamma (c,d){\Big)}\ge 
(\mu\times\sigma)(c,d).$$Indeed:\begin{eqnarray*}(\mu\times 
\sigma){\Big(}(a,b)\gamma (c,d){\Big)}&=&(\mu\times \sigma)(a\gamma 
c,b\gamma d)=\min\{\mu(a\gamma c),\sigma (b\gamma 
d)\}\\&\ge&\min\{\mu(c),\sigma (d)\} \mbox { (since } \mu,\sigma 
\mbox { are fuzzy left ideals of } M)\\&=&(\mu\times \sigma)(c,d).
\end{eqnarray*}Let now $(a,b)\le (c,d)$. Then 
\begin{eqnarray*}(\mu\times\sigma)(a,b)&=&\min\{\mu(a),\mu(b)\}\\&\le& 
\min\{\mu(c),\mu(d)\} \mbox { (since } a\le c, b\le 
d)\\&=&(\mu\times
\sigma)(c,d).\end{eqnarray*}Thus $\mu\times\sigma$ is a fuzzy left 
ideal of $M$. Similarly the cartesian product of two fuzzy right 
ideals of $M$ is a fuzzy right ideal of $M\times M$.$\hfill\Box$\\By 
Theorem 11, the following theorem holds\medskip

\noindent{\bf Theorem 12.} {\it Let M be an ordered $\Gamma$-groupoid 
and $\mu$, $\sigma$ fuzzy ideals of M. Then $\mu\times\sigma$ is a 
fuzzy ideal of $M\times M$.}\medskip

\noindent{\bf Lemma 13.} {\it Let M be an ordered $\Gamma$-groupoid 
and $\mu$, $\sigma$ fuzzy prime subsets of M. Then $\mu\times\sigma$ 
is a fuzzy prime subset of $M\times M$.}\medskip

\noindent{\bf Proof.} Let $(a,b),(c,d)\in M\times M$ and 
$\gamma\in\Gamma$. Then$$(\mu\times\sigma){\Big(}(a,b)\gamma 
(c,d){\Big)}=\max\{(\mu\times\sigma)(a,b),(\mu\times\sigma)(c,d)\}.$$In 
fact:\begin{eqnarray*}(\mu\times\sigma){\Big(}(a,b)\gamma 
(c,d){\Big)}&=&(\mu\times\sigma)(a\gamma c,b\gamma 
d)=\min\{\mu(a\gamma c),\sigma(b\gamma 
d)\}\\&=&\min{\Big\{}\max\{\mu(a),\mu(c)\}, 
\max\{\sigma(b),\sigma(d)\}{\Big\}}\\&=&\max{\Big\{}\min\{\mu(a),\sigma(b)\}, 
\min\{\mu(c),\sigma(d)\}{\Big\}}\\&=&
\max\{(\mu\times\sigma)(a,b),(\mu\times\sigma)(c,d)\}.\end{eqnarray*}$\hfill\Box$\\
By Theorem 12 and Lemma 13, we have the following theorem\medskip

\noindent{\bf Theorem 14.} {\it If M is an ordered $\Gamma$-groupoid 
and $\mu$, $\sigma$ fuzzy prime ideals of M, then $\mu\times\sigma$ 
is a fuzzy prime ideal of $M\times M$.}\medskip

\noindent{\bf Lemma 15.} {\it Let M be an ordered $\Gamma$-groupoid 
and $\mu$, $\sigma$ fuzzy semiprime subsets of M. Then 
$\mu\times\sigma$ is a fuzzy semiprime subset of $M\times 
M$.}\medskip

\noindent{\bf Proof.} Let $(a,b)\in M$ and $\gamma\in\Gamma$. Then 
$$(\mu\times\sigma)(a,b)\ge (\mu\times\sigma){\Big(}(a,b)\gamma 
(a,b){\Big)}.$$Indeed: 
$(\mu\times\sigma)(a,b)=\min\{\mu(a),\sigma(b)\}$. Since $\mu$ and 
$\sigma$ are fuzzy semiprime subsets of $M$, we have $\mu(a)\ge 
\mu(a\gamma a)$ and $\sigma(b)\ge\sigma (b\gamma b)$. Then we get
\begin{eqnarray*}(\mu\times\sigma)(a,b)&\ge&\min\{\mu(a\gamma 
a),\sigma(b\gamma b)\}=(\mu\times\sigma)(a\gamma a,b\gamma 
b)\\&=&(\mu\times\sigma){\Big(}(a,b)\gamma 
(a,b){\Big)}.\end{eqnarray*}$\hfill\Box$\\By Theorem 12 and Lemma 15, 
we have the following theorem\medskip

\noindent{\bf Theorem 16.} {\it If M is an ordered $\Gamma$-groupoid 
and $\mu$, $\sigma$ fuzzy semiprime ideals of M, then 
$\mu\times\sigma$ is a fuzzy semiprime ideal of $M\times M$.}\\By 
Theorems 4--6 and 11, 12, we have the following corollary\medskip

\noindent{\bf Corollary 17.} {\it If M is an ordered 
$\Gamma$-groupoid and $\mu$, $\sigma$ fuzzy left (resp. right) ideals 
of M, then the level subset $(\mu\times\sigma)_t$, if it is nonempty, 
is a left (resp. right) ideal of $M\times M$. If $\mu$ and $\sigma$ 
are fuzzy ideals of M and the level subset $(\mu\times\sigma)_t$ is 
nonempty, then it is an ideal of $M\times M$.}\\By Theorems 8, 10, 14 
and 16, we have the following \medskip

\noindent{\bf Corollary 18.} {\it If M is an ordered 
$\Gamma$-groupoid and $\mu$, $\sigma$ fuzzy prime (resp. semiprime) 
ideals of M, then the level subset $(\mu\times\sigma)_t$, if it is 
nonempty, is a prime (resp. semiprime) ideal of $M\times 
M$.}\medskip

As a conclusion, if $G$ is a groupoid or an ordered groupoid and 
$\mu$, $\sigma$ fuzzy left ideals, fuzzy right ideals, fuzzy ideals, 
fuzzy prime ideals or fuzzy semiprime ideals of $G$, respectively, 
then the cartesian product $\mu\times\sigma$ of $\mu$ and $\sigma$ 
is, respectively so. If $G$ is a groupoid or an ordered groupoid and 
$\mu$, $\sigma$ fuzzy left (right) ideals, fuzzy ideals, fuzzy prime 
(semiprime) ideals of $G$, respectively, then for any $t\in [0,1]$ 
for which $(\mu\times\sigma)_t$ is nonempty, the level subset 
$(\mu\times\sigma)_t$ is, respectively, a left (right) ideal, ideal, 
prime (semiprime) ideal of the cartesian product $M\times M$ of 
$M$.\bigskip

\noindent This paper has been submitted in Lobachevskii J. Math. on 
October 15, 2014.

\end{document}